\documentclass{amsart}
\usepackage{amsmath,amsfonts,amssymb,amsthm,epsfig,color}
\newcommand{\bb}{\mathbb}
\newcommand{\C}{\bb C}
\newcommand{\h}{\bb H}

\newcommand{\R}{\bb R}
\newcommand{\N}{\bb N}
\newcommand{\M}{\mathcal M}
\newcommand{\F}{\mathcal F}

\newtheorem{Theorem}[equation]{Theorem}
\newtheorem{Cor}[equation]{Corollary}
\newtheorem{Prop}[equation]{Proposition}
\newtheorem{lemma}[equation]{Lemma}
\newtheorem*{lemma*}{Lemma}
\newtheorem*{conj}{Conjecture}
\newtheorem*{claim}{Claim}
\newtheorem*{theorem*}{Theorem}

\numberwithin{equation}{section}
\begin{document}

\title[Quantitative recurrence for Teichmuller flow]{Quantitative recurrence and large deviations for Teichmuller geodesic flow}
\author[J.~S.~Athreya]{Jayadev S. Athreya\\ University of Chicago}
\subjclass[2000]{32G15, 37A10}
\email{jayadev@math.uchicago.edu}
\address{5734 S. University Avenue, Chicago, IL 60637}

\begin{abstract} 
We prove quantitative recurrence and large deviations results for the Teichmuller geodesic flow on connected components of strata of the moduli space $Q_g$ of holomorphic unit-area quadratic differentials on a compact genus $g \geq 2$ surface.

\noindent\textbf{Keywords:} Moduli spaces, geodesic flow, large
deviations.\\
\noindent\textbf{Acknowledgements:} Research partially supported by NSF grant DMS-0244542. 
\end{abstract}
\maketitle

\section{Introduction}\label{intro}

Let $\Sigma_g$ be a compact surface of genus $g \geq 2$.  Let $Q_g$ be the moduli space of unit-area holomorphic quadratic differentials on $\Sigma_g$. That is, a point $q \in Q_g$ is a equivalence class of pairs $(M, \omega)$, where $M$ is a genus $g$ Riemann surface, and $\omega$ is a holomorphic quadratic differential on $M$, i.e., a tensor with the form $f(z)dz^2$ in local coordinates, such that $\int_{M} |\omega| = 1$. Two pairs $(M_1, \omega_1)$ and $(M_2, \omega_2)$ are equivalent if there is a biholomorphism $f:M_1 \rightarrow M_2$ such that $f_{*} \omega_1 = \omega_2$.

Given a pair $q \in Q_g$, one obtains (via integration of the form) an atlas of
charts to $\C \cong \R^2$, with transition maps of the form $z \mapsto \pm z + c$. Similarly, given such an atlas of charts, one obtains a holomorphic quadratic differential by pulling back the form $dz^2$ on $\C$. 

These charts allow us to
define a $SL(2, \R)$ action on $Q_g$ (and $\tilde{Q_g}$) given by linear
post-composition with charts. 

$Q_g$ is naturally the unit cotangent bundle to $\M_g$, the moduli space of Riemann surfaces. The fiber over each point $M \in \M_g$ is the vector space of holomorphic quadratic differentials on $M$.

The space $Q_g$ is naturally stratified by integer partitions $\beta$ of $4g-4$.  Strata are not
always connected, however, they have at most finitely many
components~\cite{KZ}, and are invariant under $SL(2,\R)$. For the rest of our paper we work with one of these connected components, call it $Q$. 

Without loss of generality, we will study
strata of squares of abelian differentials (holomorphic 1-forms). Otherwise, we pass to a double cover. Since our results hold for all $g \geq 2$, we simply consider the (higher-genus) stratum defined in this way. Each stratum $Q$, while non-compact, is endowed with a continuous, ergodic, $SL(2, \R)$-invariant probability measure $\mu_Q$.

More details about the basics of moduli spaces and quadratic
differentials can be found in, e.g.,~\cite{FarbMarg, MasurTab, Strebel}.

We are particularly interested in the action of the standard subgroups $K$
(maximal compact) and $A$ (diagonal matrices) of $SL(2, \R)$. These are
both one-parameter subgroups, and can be described as follows:
$$K=\left\{r_{\theta} = \left( \begin{array}{cc} \cos \theta & \sin \theta \\
-\sin \theta & \cos \theta \end{array}\right): 0 \le \theta <
2\pi\right\}$$ 
$$A = \left\{g_t = \left(\begin{array}{cc} e^t & 0 \\ 0 & e^{-t}
\end{array}\right): t \in \R\right\}.$$We will often be interested also in the
action of the semigroup $A^{+} = \{g_t: t \geq 0\}$.

The action of $K$, known as the \emph{circle flow}, preserves the
underlying holomorphic structure, so it acts as identity when projected to $\M_g$. The action of $A$ is known as \emph{Teichmuller geodesic flow},
since the projection of an $A$-orbit yields a geodesic in the
Teichmuller metric on $\M_g$ (and in fact, all Teichmuller geodesics arise this
way). 

Masur~\cite{Masur1} and Veech~\cite{Veech1} independently showed the Teichmuller geodesic flow is ergodic with respect to $\mu_Q$ and even mixing. Veech~\cite{Veech} showed that it was `measurably Anosov', and more
recently Forni~\cite{Forni} has obtained explicit formulas regarding the hyperbolic
behavior. His results, which imply that as long as trajectories remain in
a compact set their hyperbolicity can be controlled, provide much of
the motivation for our research.

Our main results are concerned with the following scenario. Fix $q \in
Q$, and consider the
`circle' $Kq = \{r_\theta q: 0 \le \theta < 2\pi\}$. This set is endowed
with a natural probability measure $\nu$, coming from the Haar measure on $K
\cong S^1$. What is the recurrence behavior of a trajectory $\{g_t r_{\theta} q\}$ when $\theta$ is chosen at random according to $\nu$?

This type of situation was first considered by Kerckhoff-Masur-Smillie~\cite{KMS}, in order to understand dynamics of billiards in rational angle Euclidean polygons. They developed a dictionary between the dynamics of the straight-line flow on the singular Euclidean surface determined by a quadratic differential $q$ and the recurrence behavior of the geodesic trajectory $g_t q$ in the stratum $Q$. The $K$-action here corresponds to changing the direction of the straight line flow without changing the underlying surface. Thus, making a statement about a fixed $q$ and almost all directions allowed one to say something about straight line flow in almost all directions for a given flat surface, in analogy with Weyl's equidstribution theorem on the torus.

To apply these results to billiards, one follows the unfolding procedure of Zemljakov and Katok~\cite{KZ1}, which translate questions about dynamics of billiards in rational angle polygons to that of flows on an associated singular Euclidean surface. However, the set of surfaces arising from billiards is of measure 0 in every stratum, so making statements about almost every point $q$ does not suffice. To make statements about billiards, we need statements that hold for all $q \in Q$.

The main result in~\cite{KMS} is that for all $q \in Q$, and almost every $\theta \in [0, 2\pi)$, the geodesic trajectory $\{g_t r_{\theta} q\}_{t\geq 0}$ is recurrent in $Q$. As a corollary, one obtains that the directional flow for Euclidean polygonal billiards is uniquely ergodic for almost every direction. Further results in this direction include estimates on the Hausdorff dimension of divergent and bounded trajectories (\cite{Cheung, Kleinbock-Weiss}), and further relations between dynamics of the straight line flow and recurrence of the associated trajectory~\cite{Masur-Cheung1, Masur-Cheung2}.

Our main results concern finer recurrence behavior of geodesic trajectories. In particular, for any fixed $q \in Q$, we estimate the measure of the set of angles such that the associated $A^+$ orbit $\{g_t r_{\theta} q\}_{t \geq 0}$ `behaves poorly' for a length of time $T$. For us, poor behavior means that the trajectory is spending a lot of time in a neighborhood of the cusp of $Q$. Masur~\cite{Masur} proved a statistical result, known as a \emph{logarithm law} in this situation, in analogy with earlier results on symmetric spaces due to Sullivan~\cite{Sullivan} and Kleinbock-Margulis~\cite{KleinMar}.

We will construct a proper (that is, unbounded off compact sets) continuous function $V: Q \rightarrow \R^{+}$, and consider the recurrence behavior of $\{g_t r_{\theta} q\}_{t \geq 0}$ to the family of compact sets $C_l: = \{q: V(q) \le l\}$,  which form an exhaustion of $Q$ as $l$ varies. Our main results can be summarized as follows:

\begin{Theorem}\label{mainflow} Fix notation as above. Then

\begin{enumerate} 

\item For all $l$ sufficiently large and all $q \notin C_l$, there are positive constants $c_1 = c_1(l, q), c_2(l)$, with $$\nu\{\theta: g_t r_{\theta} q \notin C_l, 0 \le t \le T\} \le c_1e^{-c_2 T}$$ for all T sufficiently large. That is, the probability that a random geodesic trajectory has not visited $C_l$ by time $T$ decays exponentially in $T$.

\item For all $l, S, T$ sufficiently large and all $q \in Q$, there are positive constants $c_3 = c_3(S, l, q), c_4 = c_4(l)$, with $$\nu\{\theta: g_t r_{\theta} q \notin C_l, S \le t \le S+ T\} \le c_3e^{-c_4 T}.$$ That is, the probability that a random geodesic trajectory does not enter $C_l$ in the interval $[S, S+T]$ decays exponentially in $T$.

\item Let $q \in Q$. For any $0 < \lambda <1$, there is a
$l \geq 0$, and $0<\gamma <1$, such
that for all $T$ sufficiently large (depending on all the above constants)$$\nu\{\theta: \frac{1}{T}|\{0 \le t \le T:g_t r_{\theta} q \notin C_l\}|
> \lambda\} \le \gamma^T.$$

\end{enumerate}\end{Theorem}

Result (3) above may be thought of as a large deviations result for the Teichmuller flow. While
ergodicity guarantees that $\frac{1}{T}\left|\{0 \le t \le T: g_t q \in C_l\}\right| 
\rightarrow \mu_Q(C_l)$ for $\mu_Q$-almost every $q \in Q$, our result
gives explicit information about the likelihood of bad trajectories. Notice, however, this is \emph{not} a traditional large deviations result, which estimates the probability of a deviation of  any $\epsilon>0$ from the ergodic average.  Other interesting results concerning deviations are due to Bufetov~\cite{Bufetov}, who proved a central limit theorem for this flow.

In~\cite{Forni}, Forni related this type of fine recurrence behavior for the geodesic trajectory $\{g_t q\}_{t \geq 0}$ to deviation of ergodic averages for the straight line flow on the flat surface associated to $q$. He proved that for almost every $q \in Q$, the associated flow has the same deviation behavior. However, his result gives no information on billiards, as the set of quadratic differentials arising from billiards are a set of measure zero. Since our results provide fine recurrence information about the geodesic trajectory $\{g_t r_{\theta} q\}_{t \geq 0}$ for \emph{all} $q \in Q$, and almost all $\theta$, we conjecture the following:

\begin{conj}\nonumber For all rational-angle Euclidean polygons, the deviation of ergodic averages for the billiard flow is the same for almost all directions, and depends only on the $SL(2,\R)$ orbit of the associated quadratic differential. \end{conj}

It was also shown in~\cite{Forni} that as long as a geodesic trajectory stays within a compact set, the rate of expansion/contraction in the tangent space along the trajectory is bounded away from $1$. Thus our results can be used to obtain explicit estimates on the hyperbolicity of the flow along specific trajectories. Avila-Gouezel-Yoccoz~\cite{AGY} used exponential return estimates to a different family of compact sets to prove the exponential rate of mixing for the Teichmuller geodesic flow, which was the original motivation for this research.

Other quantitative recurrence results for dynamics on Teichmuller spaces
were obtained by Minsky-Weiss~\cite{MinskyWeiss} for the case of Teichmuller
horocycle flow.

We also have a collection of results for a certain class of random walks, defined as follows: Fix $\tau >0$.
Given that we are at a point $q \in Q$, the next point in our
trajectory will be chosen at random according to Haar measure on $S^1$
from the `circle' of radius $\tau, \{g_{\tau}r_{\theta}q: 0 \le \theta < 2\pi\}$. Note that $\mu_Q$ is a stationary measure for this walk, since it is $SL(2, \R)$-invariant. Let $\{X_n\}_{n=0}^{\infty}$ denote the random walk generated this way.

\vspace{.1in}

\noindent\textbf{Remark:} By hyperbolic geometry, one can see that trajectories of this walk
closely approximate geodesics. Thus, understanding recurrence properties
of the walk gives one insight into properties of the flow.

\vspace{.1in}

We define $P_q(E): = \mbox{Prob}(E|
X_0 = q)$ for any event $E$ defined on the trajectory $\{X_n\}_{n=0}^{\infty}$ starting at $X_0 = q$. Also, for any measurable $A \subset Q$, we define $P^n(q, A) := P_q(X_n \in A)$. That is, $P^n(q, .)$ is the probability distribution of $X_n$ given $X_0 = q$.

\begin{Theorem}\label{mainwalk} For all $\tau$ sufficiently large, we have:

\begin{enumerate}

\item There is an $l(\tau)$ such that for all $l>l(\tau)$, and all $q \notin C_l$, there are constants $c_5 = c_5(q, l,\tau), c_6 = c_6(l, \tau)$, such that $$P_q(X_i \notin C_l: 1\le i \le n) \le c_5e^{-c_6n}.$$

\item For all compact $C \subset Q$, $\epsilon >0$, there is a $C_l
\supset C$ such that $\forall q \in C$, and for all $m \geq 0$,  $$P_q(X_m
\in C_l) > 1-
\epsilon.$$

\item For all $\epsilon >0$, there is a $l = l(\epsilon)>0$ such that $\forall q \in
Q$, there is an $M(q)$ such that for all $m > M(q)$, $$P_q(X_m \in C_l)
> 1-\epsilon.$$

\end{enumerate}\end{Theorem}

The rest of the paper is organized as follows: In the next section, we give more detailed statements of our main results, and construct the function $V$. We also give a version of our results for general $SL(2, \R)$-actions. In
section~\ref{walkproof}, we prove our results for random walks. In section~\ref{polar}, we collect some technical lemmas about change of polar coordinates on hyperbolic space. In section~\ref{deviations}, we construct the required lemmas from the theory of large deviations. In section~\ref{proofs}, we prove our main theorems for the flow.

\section{Statement of Results}\label{results}

The precise statements for the flow regard a \emph{family} of proper functions $V_{\delta}$, $0 < \delta <1$, and the compact sets $C_{\delta, l} : = \{q: V_{\delta}(q) \le l\}$. We have:

\subsection{Flow results}\label{flowresults}

\begin{Theorem}\label{flowreturn} For every $1>\delta >0$ there is a proper (i.e., unbounded off compact sets), smooth, $K$-invariant function $V_{\delta}: Q \rightarrow \R^{+}$ and positive constants $t_0 = t_0(\delta), l_0 = l_0(\delta), a= a(\delta)$ such that for all $l 
\geq l_0$, there are  $1>\delta^{\prime} = \delta^{\prime}(l, \delta)>\delta$, with $\delta^{\prime}$ decreasing as a function of $l$, so that for 
all $q \notin C_{\delta, l}$ $$\nu\{\theta: g_t r_{\theta} q \notin C_{\delta,l}, 0 \le t \le
T\} \le a\frac{V_{\delta}(q)}{l}e^{-(1-\delta^{\prime}) T},$$ for all $T>t_0$.\end{Theorem}

This is the precise version of part (1) of Theorem~\ref{mainflow}. Note that the outside term essentially depends only on $V_{\delta}(q)$. 

\begin{Theorem}\label{sojourn} Let $q \in Q$. For every $\epsilon>0$, $1>\delta>0$, there are positive constants $S_0, T_1, l_1$ depending on $\delta$, $1>\delta^{\prime\prime}=\delta^{\prime\prime}(l)>\delta$, and $\alpha = \alpha(q)$ such that $$\nu\{\theta: g_t r_{\theta} q \notin C_{\delta,l}, S \le t \le S+T\} \le \alpha e^{-(1-\delta^{\prime\prime}) T},$$ for all $S> S_0$, $T>T_1$, and $l >l_1$, with $$\alpha \le 8(1+\epsilon)\sup_{\theta \in [0, 2\pi)}ab\frac{V_{\delta}(g_Sr_{\theta}q)}{l},$$ where $b$ depends only on the curvature of $\h^2$.\end{Theorem}

Here, $S_0$ depends only on the choice of curvature for the hyperbolic plane $\h^2$, and $T_1$ is the maximum of $t_0$ from Theorem~\ref{flowreturn} and a $T_2$ depending only on curvature.

\begin{Theorem}\label{proportion} Let $q \in Q$. For any $0 < \lambda <1$, and any $0 < \delta<1$ there is a
$l \geq 0$, and $\gamma <1$, such
that $$\nu\{\theta: \frac{1}{T}|\{0 \le t \le T:g_t r_{\theta} q \notin C_{\delta,l}\}|
> \lambda\} \le \gamma^T,$$ for all $T$ sufficiently
large.\end{Theorem}

This theorem uses a technical tool from the theory of large deviations, Proposition~\ref{gendev}, which makes it difficult to track the dependence of $\gamma$ on $l$ and $\lambda$. Clearly, if an $l$ works for a fixed $\lambda_0$ and $\delta$, it works for all $\lambda > \lambda_0$, and any larger $l$ (with the same $\delta$) will work for $\lambda_0$. Similarly, the same $l$ will work for the same $\lambda_0$ and any smaller $\delta$.

In a personal communication, Forni posed the following question: for \emph{every} $q$ in $Q$, and almost all $\theta$, does there exist a $\zeta>0$ such that the geodesic trajectory $\{g_t r_{\theta} q\}_{t \geq 0}$ spends at least a proportion $\zeta$ of its time in a fixed compact set? The following corollary answers in the affirmative.

\begin{Cor}\label{for} Let $q \in Q, \zeta>0$. Fix $\epsilon>0$. Let $l$ be such that Theorem~\ref{proportion} is satisfied with $\lambda = 1-\zeta, \delta = \epsilon$. Then for $\nu$-almost every $\theta$,
\begin{displaymath}\limsup_{T \rightarrow \infty} \frac{1}{T} \left|\{ 0 \le t \le T: g_t
r_\theta q \notin C_{\delta,l}\}\right| \le
\lambda.\end{displaymath}\end{Cor}

The proof of this result follows from an application of the Borel-Cantelli lemma and Theorem~\ref{proportion}.

\subsection{Random walks on $Q$}\label{walkresults}

Let $\{X_n\}_{n=0}^{\infty}$ be our trajectory. That is, $$X_{n+1} = g_{\tau} r_{\theta_{n}} X_n,$$ where $\{\theta_n\}_{n=0}^{\infty}$ are independent and identically distributed (i.i.d.) according to the uniform distribution on $[0, 2\pi)$ (equivalently, $X_{n+1} = g_{\tau} k_n X_n$, where $\{k_n\}_{n=0}^{\infty}$ are i.i.d. according to Haar measure on $K$). This implies $\{X_n\}$ is a Markov chain, with transition probability function $P(x, A) = \nu\{\theta: g_{\tau}r_{\theta}q \in A\}$. 

\begin{Theorem}\label{walkreturn} Fix $0<\delta<1$. Let $\{C_{\delta, l}\}$ be as in
Theorem~\ref{flowreturn}. Then there is a $\tau_0 = \tau_0(\delta)>0$ such that for all $\tau > \tau_0$, there is a $\tilde{l}_0 = \tilde{l_0}(\delta, \tau)$ such that for all $l >\tilde{l_0}$, there is a $\gamma = \gamma(l, \delta)<1$ such that
for all $q \notin C_{\delta,l}$, and for all $n \geq 1$, $$P_q\left(X_j \notin C_{\delta,l}: 0 \le j \le n\right) \le \frac{V_{\delta}(q)}{l}\gamma^n.$$ \end{Theorem}

This is completely analagous to Theorem~\ref{flowreturn}, and indeed the proof of
this result will be essential to the proof of Theorem~\ref{flowreturn}.

For the rest of the section, fix $1>\delta>0$, and $\tau > \tau_0(\delta)$.

\begin{Theorem}\label{tightness}
Let $\epsilon >0$. Then:
\begin{enumerate}
\item For all compact $C \subset Q$, there is a $C_{\delta,l}
\supset C$ such that $\forall q \in C$, and for all $m \geq 0$,  $$P_q(X_m
\in C_{\delta, l}) > 1-
\epsilon.$$
\item There is a $l>0$ such that $\forall q \in
Q$, there is an $M(q)$ such that for all $m > M(q)$, $$P_q(X_m \in C_{\delta,l})
> 1-\epsilon.$$
\end{enumerate}
\end{Theorem}

%We also have a mixing results for these walks. For an arbitrary signed measure $\upsilon$ on $Q$, we define $$||\upsilon||_V := \sup_{g: |g| \le V} |\upsilon(g)|,$$ where the supremum is taken over all measurable functions.

%\begin{Theorem}\label{walkmixing} There are positive constants $c_{7}>1, c_{8}$ such that for almost every $q \in Q$,  $$\sum_{n} c_{7}^n ||P^n(q, .) - \mu_Q ||_{V_{\delta}} \le c_{8}V_{\delta}(q).$$\end{Theorem}

%This is a strong statement about rapid convergence of these random walks to the stationary distribution. Note that $\mu_Q$ is \emph{not} the unique stationary distribution for these walks. Any $SL(2, \R)$ invariant measure on $Q$ is invariant under the walk.

Finally, we have:

\begin{Cor}\label{walkfor}
\item Fix $q \in Q$, $0<\lambda <1$, and $\epsilon >0$. Then there is
a $ l=l(q, \lambda, \delta, \epsilon)>0$ so that, for all $n$ sufficiently large,
\begin{displaymath} P_q\left(\frac{1}{n}|\{1 \le i \le n: X_i \in
C_{\delta,l}\}| > \lambda \right) >1-\epsilon \end{displaymath}\end{Cor}
 
This is analogous to Corollary~\ref{for}.
 
\vspace{.1in} 
 
\noindent \textbf{Remark:} $l$ can be chosen uniformly as $q$ varies
over a compact set.   

\vspace{.1in}

\subsection{$SL(2, \R)$-actions}\label{axiomatic}

Let $SL(2, \R)$ act continuously on a topological space $X$. Suppose there is a family of $K$-invariant functions $V_{\delta}: X \rightarrow \R$, $0 < \delta <1$ satisfying the following properties:

\begin{itemize}

\item  For all $x \in X$, consider the function $V_{\delta, x}: SL(2,
\R) \rightarrow \R^+$ defined by $V_{\delta,x} = V_{\delta}(gx)$. Note that by $K$-invariance, we can view $V_{\delta,x}$ as a function on $\h^2 = SO(2)$\verb"\"$SL(2, \R)$. We require the following : For all $\sigma >1$, there exists a $\kappa>0$ such that for any $p \in \h^2$
with $d(p, i) < \kappa$,
and any $h \in SL(2, \R)$\begin{equation}\label{logsmooth} \sigma^{-1}V_{\delta,x}(h)
\le V_{\delta,x}(ph) \le \sigma V_{\delta,x}(h).\end{equation} We say such a function is \emph{logsmooth}.

\item 
For all $1>\delta>0$, there is a constant $\tilde{c} = \tilde{c}(\delta)$ such that for all sufficiently large $\tau$, there is a $\tilde{b} = \tilde{b}(\tau,\delta)$
such that  for all $x \in X$,  \begin{equation}\label{drift1} (A_\tau V_{\delta})(x) : = \int_{0}^{2\pi}
V_{\delta}(g_{\tau} r_{\theta} x) d\nu(\theta)  \le \tilde{c}e^{-(1-\delta)\tau}V_{\delta}(x) + \tilde{b}.\end{equation}

\end{itemize}

All our results on recurrence to the sublevel sets $C_{\delta, l}$ from subsections~\ref{flowresults} and~\ref{walkresults}, while constructed for the space $Q$ hold for any space $X$ with the above properties.

The key lemma needed to prove our results will be the construction of a family of functions $\{V_{\delta}\}$ on $Q$ satisfying the above requirements.

\begin{lemma}\label{drift}There is a family of smooth, proper functions
$V_{\delta}: Q \rightarrow \R^{+}$, $0 < \delta <1$ satisfying equations ~\eqref{logsmooth} and ~\eqref{drift1}.
\end{lemma}

We require the following technical lemma from~\cite{EskinMasur} (page
465, Lemma 7.5). 

\begin{lemma*}For all $\delta>0$ there are logsmooth functions $V_0, \ldots ,V_n:
Q \rightarrow \R^+$ such that $V_0$ is proper, and for every $\tau>0$, there are constants $w= w(\tau, \delta),
\tilde{b^{\prime}}=
\tilde{b^{\prime}}(\tau, \delta)$, and $\tilde{c^{\prime}} = \tilde{c^{\prime}}(\delta)$, \emph{independent} of $\tau$
such that for all $0 \le i \le n$, and $\forall q \in Q$,
\begin{displaymath} (A_{\tau}V_i)(q) \le
\tilde{c^{\prime}}e^{-(1-\delta)\tau}V_i(q) + w\sum_{j= i+1}^n V_j(q) +
\tilde{b^{\prime}}.\end{displaymath}\end{lemma*}

\vspace{.1in}

\noindent\textbf{Remark:} In fact $V_0(q) = \max\left(1, \frac{1}{l(q)^{1+\delta}}\right)$, where $l(q)$ denotes the length of the shortest saddle connection on $q$. Recall that a \emph{saddle connection} is a geodesic (in the metric determined by $q$) connecting two zeroes of $q$.

\vspace{.1in}

\noindent\textbf{Proof of Lemma~\ref{drift}:}

\noindent Fix $\delta, \tau>0$. Let $\tilde{c} = 2\tilde{c^{\prime}}$
$\lambda_0 = \frac{w}{\tilde{c^{\prime}}}$, $\lambda_i = \left(\frac{\tilde{c^{\prime}}}{w}+1\right)^{i-1}$, for $0 \le i \le n$. Note that \begin{displaymath}
\sum_{i=0}^{j-1} \lambda_i \le
\frac{\tilde{c}}{2w}\lambda_j.\end{displaymath}\noindent Set $\tilde{b} =
\tilde{b^{\prime}}\sum_{i=0}^n\lambda_i$. Let $V_{\delta}(q) =
\sum_{i=0}^{n} \lambda_i V_i(q)$. Then 

\begin{eqnarray}
\left(A_\tau V_{\delta}\right)(q) = \sum_{i=0}^n \lambda_i
\left(A_{\tau}V_i\right)(q)\nonumber &\le& \sum_{i=0}^n
\lambda_i\left(\tilde{c^{\prime}}e^{-(1-\delta)\tau}V_i(q) + w\sum_{j=i+1}^n V_j(q) +
\tilde{b^{\prime}}\right) \nonumber \\ & =& \sum_{i=0}^n \lambda_i
\tilde{c^{\prime}} e^{-(1-\delta)\tau}
V_i(q) + w\sum_{j=1}^n \left(\sum_{i=1}^{j-1}
\lambda_i\right)V_j(q)+ b \nonumber \\ &=&
\tilde{c^{\prime}}e^{-(1-\delta)\tau}\lambda_0V_0(q) + \sum_{i=1}^n
V_j(q)\left(\lambda_j \tilde{c^{\prime}}e^{-(1-\delta)\tau} +
w\sum_{i=0}^{j-1}\lambda_i\right) + \tilde{b} \nonumber \\ &=&
\tilde{c^{\prime}}e^{-(1-\delta)\tau} \lambda_0 V_0(q) + \sum_{j=1}^n
2\lambda_j\tilde{c^{\prime}}e^{-(1-\delta)\tau}V_j(q) + \tilde{b} \nonumber \\ &\le&
\tilde{c}e^{-(1-\delta)\tau}V_{\delta}(q) + \tilde{b}. \nonumber\end{eqnarray} 

Thus we have constructed a family of functions $V_{\delta}$ satisfying equations~\eqref{logsmooth} and~\eqref{drift1}. That $V_{\delta}$ satisfies~\eqref{logsmooth} follows from the logsmoothness of the $V_i$'s  noted in~\cite{EskinMasur} (where they are called $\alpha_i$'s) on page 471, at the beginning of the proof of Proposition 7.2. Finally, $V_{\delta}$ is proper since $V_0$ is proper. \qed\medskip

\noindent 

\section{General Markov Chain results}\label{walkproof} 

Our goal in this section is to recall some results from the theory of Markov Chains which will allow
us to prove Theorems~\ref{walkreturn}-\ref{tightness}. We also give some applications to
random walks on homogeneous spaces, as considered in~\cite{EskinMarg}. Many of
these results can be found, in greater generality, in~\cite{MeynTweedie}.

Before stating our main results, we recall some basic notation: If
$\{X_n\}$ is a Markov chain on $(S, \mathcal{S})$, $S$ the state space, and
$\mathcal{S}$ the $\sigma$-algebra, then, for any event $E$ (an event $E$ is a
set in the product $\sigma$-algebra) and any starting point $x \in S$ ,
$P_x(E) : = P(E| X_0 = x)$, i.e., it is the probability of the event $E$
occuring given that our starting point was $x$. Similarly, given a measurable subset $C \in S$, we write $P^n(x, C) = P_x(X_n \in C)$, and we define the first hitting time of $C$ by $\tau_C := \inf\{n \geq 1: X_n
\in C\}$. 

%We make the following definition:

%\begin{Def} $\{X_n\}$ satisfies the \emph{minorization condition} with respect to the stationary measure $\pi$ on $(S, \mathcal{S})$ if there is a probability measure $\rho$ so that:

%\begin{itemize}

%\item There is a set $A$, with $\rho(A) =1$, and an $\epsilon>0$ such that for all $B \in \mathcal{S}$, $$\inf_{x \in A} P(x, B) > \epsilon\rho(B).$$

%\item For $\pi$-almost all $x \in S$, $$\pi\{x : P_x(\tau_A < \infty)>0\} = 1.$$

%\end{itemize}
%\end{Def}

The main result of this section is:

\begin{Prop}~\label{walkgen}Let $S$ be a non-compact topological space, and $\mathcal{S}$ its
Borel $\sigma$-algebra. Let $\{X_n\}_{n=0}^{\infty}$ be a Markov chain on
$(S, \mathcal{S})$. Suppose there exists a smooth, proper function $V: S
\rightarrow \R^{+}$ and constants $0<c<1$ and $b>0$ such that $$(PV)(x) :
= E(V(X_1)|X_0 =x) \le cV(x) + b,$$ for all $x \in S$. Then
\begin{enumerate}

\item For any $l >0$ and $x \notin C_{l}: = \{y
\in S: V(y) \le l\}$, $$p_n(x) := P_x\left(\tau_{C_{l}} > n\right) \le
\frac{V(x)}{l} \left(c+\frac{b}{l}\right)^n,$$ for all $n \geq 0$.

\item For all $\epsilon >0$, and for all compact $C \subset S$, there is
an $l$ such that $\forall x \in C, m \geq 0$, $$P^m\left(x, C_l\right) >
1-\epsilon.$$

\item For all $\epsilon >0$ there is an $l >0$ such that $\forall x \in
S$, there is an $M(x) >0$ such that for $m> M(x)$, $$P^m\left(x, C_l\right) >
1- \epsilon.$$

%\item If in addition $\{X_n\}$ satisfies the minorization condition with respect to the measure $\pi$, there are positive constants $r>1, R>0$ such that for $\pi$-almost every $x \in S$, $$\sum_n r^n ||P^n\left(x, .\right) - \pi||_V \le RV(x),$$where $||\upsilon||_V := \sup_{g : |g| \le V} |\upsilon(g)|$ (supremum taken over $\upsilon$-measurable $g$) for any signed measure $\upsilon$.

\item  For all $x \in S$, $0<\lambda <1$, and $\epsilon>0$, there is
a $ l=l(q, \lambda, \epsilon)>0$ so that, for all $n$ sufficiently large,
\begin{displaymath} P_x\left(\frac{1}{n}|\{1 \le i \le n: X_i \in
C_l\}| > \lambda \right) >1-\epsilon .\end{displaymath} If we fix $\lambda$ and $\epsilon$, then $l$ can be chosen uniformly as $x$ varies over a compact set.

\end{enumerate}

\end{Prop}

\noindent\textbf{Proof of Theorems~\ref{walkreturn}-\ref{tightness} and Corollary~\ref{walkfor} :} Combine Proposition~\ref{walkgen} with Lemma~\ref{drift}. For Theorem~\ref{walkreturn}, fix $\delta >0$ and let $\tau_0$ be such that $c = \tilde{c}e^{-(1-\delta)\tau_0} <1$, we obtain our result with $\gamma = \left(c + \frac{\tilde{b}}{l}\right)$ for $l > \frac{\tilde{b}}{1-c}$.\qed\medskip

\noindent\textbf{Proof of Proposition~\ref{walkgen}:} 

\begin{itemize}

\item\textbf{Proof of (1):} Let $B_n : = \{\tau_{C_l} > n\}$. Then $p_n(x) = P_x(B_n)$. For
$n \geq 0$, \begin{displaymath}l p_n \le E_x(V(X_n): B_n) =: D_n
\end{displaymath} since on $B_n$, $X_n \notin C_l$, i.e.
$V(X_n) \geq l$.

Now, $B_n \subset B_{n-1}$, so 
\begin{eqnarray} 
D_n &\le& E_x(V(X_n): B_{n-1}) \nonumber\\ 
&=& E_x(E(V(X_n)|X_{n-1}): B_{n-1}) \\
&=& E_x( (PV)(X_{n-1}) : B_{n-1}) \nonumber \label{markov}
\end{eqnarray} where we are using the Markov property in the 2nd line.

Now we can apply our 
condition $$(PV)(X_{n-1}) \le cV(X_{n-1}) + b.$$ This yields $$D_n \le cD_{n-1} +
bp_{n-1},$$and using the observation that $p_{n-1} \le \frac{D_{n-1}}{l}$, we obtain the recurrence relation $$D_n \le \left(c+\frac{b}{l}\right)D_{n-1}.$$Iterating this, we obtain $$D_n \le D_0 \left(c+ \frac{b}{l}\right)^n.$$ Since $D_0 = V(X_0) = V(x)$, and $p_n \le \frac{D_n}{l}$, we obtain our result. Note that the result is only meaningful if $(c+\frac{b}{l}) < 1$, which is equivalent to setting $l > \frac{b}{1-c}$.

\item\textbf{Proof of (2):} We have $(PV)(x) \le cV(X) + b$ for $c<1$.
Iterating this,
we get that $(P^mV)(x) \le c^mV(x) + b^{\prime}$, where $b^{\prime}$ does
not depend on $m$ or $x$. Set $l = \sup_{y \in C}\frac{V(y)+b^{\prime}}{\epsilon}$.
Then we have that \begin{displaymath} lP^m(x, C_l^c) \le E_x(V(X_m)) \le c^mV(x) + b^{\prime},\end{displaymath} so
we get that $P^m(x, C_l^c) < \epsilon$ as desired. \qed\medskip

\item\textbf{Proof of (3):} For the third property, we select l = $2b^{\prime}/\epsilon$, where $b^{\prime}$ is as above. For $m$ sufficiently large, $c^mV(x) \le
b^{\prime}$, so by the argument in part (2), we can get our conclusion. \qed\medskip

%\item\textbf{Proof of (4):} Combine ~\cite{MeynTweedie}, Theorem 15.0.1, page 359 with the results of ~\cite{ADS}, we obtain our result.\qed\medskip

\item\textbf{Proof of (4):} Fix $x \in S$,
$0<\lambda<1$, and $\epsilon>0$. By part $(3)$, there is
an $l>0$ so that for $n$ sufficiently large, $P(X_n \in C_l) > 1-
\epsilon^{\prime}$, where  $0<\epsilon^{\prime} <
\frac{2}{3}\epsilon(1-\lambda)$. Set $S_n = \frac{1}{n}\sum_{i=1}^n
\chi_{C_l}(X_i)$, where $\chi_{C_l}$ is the indicator function of $C_l$.
Then, for any $\lambda <1$, 
\begin{eqnarray} E(S_n) &\le& \lambda P(S_n
\le \lambda) + P(S_n > \lambda) \nonumber\\ &=& \lambda +
(1-\lambda)P(S_n > \lambda). \nonumber
\end{eqnarray}

Thus, we have \begin{displaymath} P(S_n>\lambda) \geq \frac{E(S_n) -
\lambda}{1-\lambda}. \end{displaymath} 

Now, for $n$ sufficiently large $E(S_n) \geq 1-\frac{3}{2}\epsilon^{\prime}$, thus,
\begin{displaymath}P(S_n
> \lambda) \geq \frac{1-\frac{3}{2}\epsilon^{\prime} -\lambda}{1-\lambda} >
1-\epsilon.\end{displaymath}

The fact that $l$ can be chosen uniformly as $x$ varies over a compact
set follows from part (2).\qed\medskip

\end{itemize}

These types of questions were considered for random walks on homogeneous spaces in~\cite{EskinMarg} by Eskin and Margulis. They constructed a function $V$ on their state space satisfying the conditions of Proposition~\ref{walkgen}, and used this to draw conclusions (2) and (3). Conclusions (1) and (4) appear to be new results for these walks.

Precisely, we have the following:

\begin{theorem*} Let $G$ be a semisimple Lie group, and $\Gamma$ a non-uniform lattice. Let $\mu$ be a probability measure on $G$ satisfying the conditions of Theorem 2.1 in~\cite{EskinMarg}. Consider the Markov chain $\{X_n\}_{n=0}^{\infty}$ defined on $G/\Gamma$ by the measure $\mu$: $$X_{n+1} = g_n X_n,$$ with $\{g_n\}_{n=0}^{\infty}$ an i.i.d. (with distribution $\mu$) sequence of elements of $G$. Then there is a function $V: G/\Gamma \rightarrow \R^{+}$ satisfying the conditions of Proposition~\ref{walkgen}. Thus, conclusions (1)-(4) of the Proposition are satisfied.\end{theorem*}

\section{Polar coordinates and shadowing}\label{polar}

We require two lemmas about change of polar coordinates in the hyperbolic plane $\h^2 =SO(2)$\verb"\"$SL(2, \R)$. We fix two positive numbers $t_1, t_2$, and basepoints $i$ and $z_0 = i. g_{t_1} r_{\theta}$ (these will correspond to our basepoint $q$ and an arbitrary $q_0$ in its $SL(2, \R)$-orbit, projected to $SO(2)$\verb"\"$Q$). We let $d(., .)$ denote distance in the hyperbolic plane.

Consider the circle of radius $t_2$ around $z_0$, defined by $\{z_{\phi} = i.g_{t_2}r_{\phi}g_{t_1} r_{\theta}: 0 \le \phi <2\pi\}$. We say that $t_2, \phi$ are the polar coordinates of $z_{\phi}$ based at $z_0$. For each $\phi$, we define $D = D_{t_1, t_2}(\phi)$ and $\Psi = \Psi_{t_1,t_2}(\phi)$ by $z_{\phi} = i. g_{D(\phi)} r_{\theta+\Psi(\phi)}$, i.e. $D(\phi), \theta+\Psi(\phi)$ are the polar coordinates of $z_{\phi}$ based at $i$. Note that $D, \Psi$ are \emph{independent} of $\theta$.

Geometrically, $D(\phi)$ is the
distance and $\Psi(\phi)$ is the angle (measured clockwise from the the geodesic connecting $i$ to $z_0$) of the
geodesic segment connecting $i$ to $z_{\phi}$. Hyperbolic trigonometry (the laws of sines and
cosines, appied to the triangle formed by the points $i$, $z_0$, and
$z_{\phi}$) yield: \begin{equation}\label{polarlength} \cosh D(\phi) =
\cosh t_1 \cosh t_2+ \sinh t_1\sinh t_2 \cos \phi ,\end{equation} and
\begin{equation}\label{polarangle} \sin \Psi(\phi) = \frac{\sinh t_2}{\sinh
D(\phi)} \sin \phi.\end{equation}

If $t_2>t_1$, $i$ lies inside the circle of radius $\tau$ around $z_0$, and
thus, the map $\Psi$ is both one-to-one and onto. If $t_1>t_2$, the point $i$ is
outside the cirlce, and $\Psi$ is neither one-to-one or onto. In this case,
the image is an interval, with boundary points such that the geodesic
determined by those angles intersects the circle of radius $t$
tangentially. In the interior of the interval, each point $\psi$ has two
preimages, call them $\phi_1, \phi_2$ one such that $D(\phi_1) \approx
t_1-t_2$, and one such that $D(\phi_2) \approx t_1+t_2$. For our applications, we
will only be concerned with $\phi_2$. 

The key technical lemma is as follows:

\begin{lemma}\label{shadow}
Let $A \subset [0, 2\pi)$ be a measurable set. Then, for every $\epsilon>0$, there are $\tau_1, \tau_2$, such that for all $t_1 > \tau_1, t_2> \tau_2$ the neighborhood $U = \Psi_{t_1, t_2}([-\pi/2, \pi/2])$ of $0$ satisfies \begin{displaymath} \frac{\nu(U \cap
\Psi(A))}{\nu(U)} \le 4(1+\epsilon)\nu(A).\end{displaymath}
\end{lemma}

\noindent\textbf{Proof:} For this estimate, we need to control the behavior of the
derivative $\Psi^{\prime}$. More precisely, we need to control \emph{ratios}
$\Psi^{\prime}(\phi_1)/\Psi^{\prime}(\phi_2)$, with $\phi_1, \phi_2 \in \Psi^{-1}
U$, so we can compare $\nu$ and $\Psi_{*} \nu$, where $\Psi_{*}\nu(E) = \nu(\Psi^{-1}E)$. 

We have the following claim:

\begin{claim}\nonumber Let $\eta >0$. For $t_1, t_2$ sufficiently large, $$\frac{e^{-t_1}}{2}(1-\eta) \le |\Psi^{\prime}(\phi)| \le e^{-t_1}(1+\eta),$$ for $\phi \in [-\pi/2, \pi/2]$.\end{claim}

\noindent \textbf{Proof of Claim:}

Implicit differentiation of equations~\eqref{polarlength} and~\eqref{polarangle}
yield:

\begin{equation}\label{lengthderiv} D^{\prime}(\phi) \sinh D(\phi) =
-\sinh t_1 \sinh t_2 \sin \phi \end{equation}

and 

\begin{equation}\label{anglederiv} \Psi^{\prime}(\phi)\cos \Psi(\phi) = \sinh t_2
\frac{\cos \phi \sinh D(\phi) + \sin^2 \phi \coth D(\phi) \sinh t_1 \sinh
t_2}{\sinh^2 D(\phi)}.\end{equation}

Let $\kappa >0$. Let $t_1, t_2$ be large enough so that
 for all $\phi \in [-\pi/2, \pi/2]$, 
 \begin{enumerate}
\item $D(\phi) > t_1 + t_2 - \kappa$, 
\item $\coth(D(\phi)) > 1-\kappa$
\item $\cos(\Psi(\phi))>1-\kappa$
\item $|\sinh t_1 \sinh t_2 - \sinh D(\phi)| \le \kappa$
\item $1-\kappa<\frac{2\sinh t_1}{e^{t_i/2}}, \frac{2\cosh t_1}{e^{t_i/2}}<1+\kappa$ for $i=1, 2$.
\end{enumerate}

That we can achieve the above inequalities follows from hyperbolic geometry and the basic properties of $\sinh$ and $\cosh$.

Let $\eta >0$. Using the above inequalities and some basic algebra, we can select $\tau_1, \tau_2$ such that for all $t_1>\tau_1, t_2>\tau_2$, we have

$$(1-\eta) \le \frac{\Psi^{\prime}(\phi)}{ \frac{e^{-t_1} }{2}(2\cos\phi +\sin^2 \phi)} \le (1+\eta).$$

\vspace{.1in}

\noindent \textbf{Remark:} The expression $\frac{e^{-t_1} }{2}(2\cos\phi +\sin^2 \phi)$ is obtained by replacing the quantities in equation~\eqref{anglederiv} with their approximations (1)-(5). 

\vspace{.1in}

Now, let $f(\phi) = 2\cos\phi + \sin^2 \phi$. We have $1 \le f(\phi) \le 2$ for $\phi \in [-\pi/2, \pi/2]$. Thus, $$\frac{e^{-t_1}}{2}(1-\eta) \le \Psi^{\prime}(\phi) \le (1+\eta)e^{-t_1},$$ completing the proof of the claim. 

To complete the proof of the lemma, let $\epsilon >0$ and $\eta$ be such that $\frac{1+\eta}{1-\eta} \le 1+ \epsilon$. By the claim, we know
that the proportion of
measure of any set in $[-\pi/2, \pi/2]$ cannot be expanded by more than
$2\frac{1+\eta}{1-\eta} \le 2(1+\epsilon)$ under $\Psi$, since that is the maximum of $\frac{\Psi^{\prime}(\phi_1)}{\Psi^{\prime}(\phi_2)}$ for $\phi_1, \phi_2 \in [-\pi/2, \pi/2]$. Now, since $\nu([-\pi/2, \pi/2]) = 1/2$, we have
\begin{displaymath} \frac{\nu(A \cap [-\pi/2, \pi/2])}{\nu([-\pi/2, \pi/2])}
\le 2\nu(A), \end{displaymath} which yields \begin{displaymath}
\frac{\nu(\Psi(A) \cap U)}{\nu(U)} \le 2(1+\epsilon) \cdot 2\nu(A) = 4(1+\epsilon)\nu(A).\end{displaymath}\qed\medskip

Our second main lemma is as follows: let $\tau >0$. For any $t>0$, define $$I(t) : = \{0 \le \theta < 2\pi: d(i.g_t, i.g_t r_{\theta}) \le 3\tau\}.$$ Fix $\kappa>0$. For any $\theta \in I(t)$, define $$U_{\theta} := \{ 0 \le \phi < 2\pi: D_{t, \tau}(\phi) > 1-\kappa, (\Psi_{t, \tau}(\phi) + \theta) \in I(t+\tau)\}.$$ Let $L_{\theta} = \nu(U_{\theta})$. The following lemma is proved in~\cite{EskinMasur}. For notational convenience, we write drop the subscripts for $T$ and $S$.

\begin{lemma}\label{polar2} There is a constant $c^{\prime\prime}>0$ such that for all $\kappa>0$, there is a $\tau>0$ such that for all $\theta \in I(t)$, the map $\Psi|_{L_{\theta}}$ is a diffeomorphism onto its image, and, making the subsitution $\psi = \Psi(\phi)$, we have $$c^{\prime\prime} \le \nu(L_{\theta}) \le \int_{L_{\theta}} d\phi  = \int_{\Psi(L_{\theta})} |\frac{d\phi}{d\psi}| d\psi.$$\end{lemma}

\noindent\textbf{Proof:}\cite{EskinMasur}, page 467, Lemma 7.6.

\vspace{.1in}

\noindent\textbf{Remarks:}
\begin{itemize}
\item $d\theta$ denotes $d\nu(\theta)$, and since we have normalized $\nu$ to be a probability measure, we do not need to divide by $2\pi$.
\item In~\cite{EskinMasur}, the sets $I(t)$ are defined by $I(t): = [-\rho e^{-t}, \rho e^{-t}]$ for some positive constant $\rho$, and they require $\rho$ to be large enough so that the diameter of the set $J(t) = \{i.g_t r_{\theta}: \theta \in I(t)\}$ is at least $2\tau$. By hyperbolic geometry, our sets $I(t)$ are of this form, and they obviously satisfy the required condition.
\end{itemize}

\vspace{.1in}

We have the following Corollary, also from~\cite{EskinMasur}: 

\begin{Cor}\label{polar3} Let $f:SL(2, \R) \rightarrow \R$ be a logsmooth $SO(2)$-invariant function. Fix $\sigma > 1$. Let $\kappa>0$ be as in equation~\eqref{logsmooth}. Fix $\tau$ so that Lemma~\ref{polar2} holds. Then there is a $c^{\prime} >0$, independent of $\tau$, such that $$\int_{I(t+\tau)} f(g_{t+\tau}r_{\theta})d\theta \le c^{\prime}\sigma\int_{I(t)} (A_{\tau}f)(g_t r_{\theta})d\theta.$$\end{Cor}

\noindent\textbf{Proof:}\cite{EskinMasur}, page 468, Lemma 7.7.

\section{Large Deviations}\label{deviations}

In this section, we prove the key technical lemma for our main large
deviations result Theorem~\ref{proportion}. We assume some familiarity
with the theory of conditional expectation. Excellent references
include~\cite{Varadhan, Durrett}. 

\begin{Prop}\label{gendev} Let $\{\tau_i\}_{i=0}^{\infty}$ be a sequence of positive real-valued random
variables on a probability space $(\Omega, \F, P)$. Let $\{\F_i\}$ be af
filtration of $\F$ such that for all $i$, $\tau_i \in \F_i$, i.e.,
$\tau_i$ is $\F_i$-measurable. Suppose there exist positive random
variables $\eta, \xi$ with $E\eta < E\xi$, $E\xi>0$, and a real number $\theta_0>0$
such that, for all $0 \le \theta < \theta_0$:
\begin{enumerate}
\item $E(e^{\theta \tau_{2i}}|\F_{2i-1}) \le E(e^{\theta \eta})$
\item $E(e^{-\theta (\tau_{2i-1}+\tau_{2i})}|\F_{2i-1}) \le E(e^{-\theta \xi})$
\end{enumerate}

Let $T_n = \sum_{i=0}^n \tau_i$. Let \begin{equation}X(t) = \left\{
\begin{array}{ll}  1 & T_{2i-1} \le t
< T_{2i}\\ 0 & \mbox{otherwise}\end{array}\right.\nonumber\end{equation}

Then, $\forall \lambda > E\eta/E\xi$, there is a $0 < \gamma < 1$ such
that for all $T$ sufficiently large\begin{displaymath} P\left(\frac{1}{T}\int_0^T X(t)dt >
\lambda\right) \le
\gamma^T.\end{displaymath}
\end{Prop}

\vspace{.1in}

\noindent For our application, $\tau_{2i}$ will be the time a random geodesic spends
outside the compact set, and $\tau_{2i-1}$ the time inside on the $i$th
`sojourn'. $\eta$ is a stochastic upper bound for $\tau_{2i}$, and $\xi$ as a stochastic lower
bound for the length of a `cycle'. Thus condition 1 should be thought of as a stochastic upper bound for the time spent outside and condition 2 a stochastic lower bound for the total time of a
cycle. In our application, we will give a stronger, deterministic lower
bound: in fact, we will construct our variables so that $\tau_{2i-1} > C$, for some fixed $C$. This clearly implies condition 2, simply by taking $\xi = C$. For condition 1, we will use Theorem~\ref{flowreturn} to show that $\tau_{2i}$ cannot grow too large.

\vspace{.1in}

\noindent
\textbf{Proof:}

Define $f(\theta):=E(e^{\theta \eta})$ and $g(\theta):= E(e^{-\theta \xi})$.

Let $N(T) = \sup\nu\{k: T_{2k} \le T\}$. Then,
\begin{eqnarray}
P\left(\int_0^T X(t)dt>\lambda T\right)\nonumber &\le& P\left(\int_0^T
X(t)dt > \lambda T \mbox{ and } N(T) \le cT\right) + P\left(N(T) >
cT\right)\nonumber \\
&=& P\left(\sum_{i=0}^{\lfloor cT/2 \rfloor} \tau_{2i} > \lambda T\right) +
P\left(N(T) > cT\right).\end{eqnarray}

We estimate each of these terms in turn. Let
$n:=\lfloor cT/2 \rfloor$. Then, since $$P\left(\sum_{i=0}^{\lfloor cT/2 \rfloor} \tau_{2i} > \lambda^{\prime}n\right) \geq P\left(\sum_{i=0}^{\lfloor cT/2 \rfloor} \tau_{2i} > \lambda T\right) \geq  P\left(\sum_{i=0}^{\lfloor cT/2 \rfloor} \tau_{2i} > \lambda^{\prime}(n+1)\right)$$ to estimate the first term it suffices to estimate
$P\left(\sum_{i=0}^{n} \tau_{2i}  > \lambda^{\prime} n\right)$,
where $\lambda^{\prime} =2\lambda/c$. Now, for any $\theta_ 0 > \theta
\geq 0$ 
\begin{eqnarray} 
P\left(\sum_{i=0}^n \tau_{2i} > \lambda^{\prime}n\right) &=& P\left(e^{\theta
\sum_{i=0}^n\tau_{2i}} > e^{\theta\lambda^{\prime}n}\right) \nonumber\\
&\le& e^{-\theta \lambda^{\prime} n} E\left(e^{\theta \sum_{i=0}^n
\tau_{2i}}\right)
\nonumber\\
&\le& e^{-\theta \lambda^{\prime} n} f(\theta)^n. \label{cond}
\end{eqnarray}

The last inequality follows from the first condition in our theorem, and
the fact that each $\tau_i$ is $\F_i$-measurable. Since
equation~\label{expec} holds for any $\theta_0 > \theta \geq 0$, we have
\begin{displaymath} P\left(\sum_{i=0}^n \tau_{2i} >
\lambda^{\prime}n\right) \le
\inf_{\theta_0 > \theta \geq 0} \left(f(\theta)e^{-\theta
\lambda^{\prime}}\right)^n.\end{displaymath} 

Let $\lambda^{\prime} > E\eta$. Then, letting \begin{displaymath}F(\theta) =
f(\theta) e^{-\theta\lambda^{\prime}} = E\left(e^{\theta(\eta
-\lambda^{\prime})}\right),\end{displaymath} we get $F(0) = 1$,
\begin{displaymath} F^{\prime}(\theta) = E\left(e^{-\theta
(\eta - \lambda^{\prime})} (\eta - \lambda^{\prime})\right). \end{displaymath}
This implies that \begin{displaymath} F^{\prime}(0^+) = E\left(\eta
-\lambda^{\prime}\right) < 0.\end{displaymath}

\noindent Thus, there is a $0 \le \theta_1 \le \theta_0$ such that
$F(\theta_1) : = \gamma^{\prime} < 1$. Plugging this into equation~\eqref{cond} yields the estimate for our first term. 

To estimate the second term, let $\xi_i = \tau_{2i-1} + \tau_{2i}$. Fix
$c > 1/E\xi$. By a similar argument to that above, we obtain
\begin{displaymath} P\left(N(T) > cT\right) \le \inf_{0 \le \theta \le \theta_0}
e^{\theta T} g(\theta)^{cT}.\end{displaymath}

\noindent Let \begin{displaymath}G(\theta) = g(\theta)e^{\theta/c} =
E\left(e^{-\theta(\xi -1/c)}\right).\end{displaymath} Once again, as above, we
obtain \begin{displaymath} G^{\prime}(0^+) = -E\left(\xi - 1/c\right) \le
0.\end{displaymath} Thus, there exists $\theta_2$ with $G(\theta_2) =
\gamma^{\prime\prime} <1$, and so we have our desired estimate.\qed\medskip

\begin{Cor}\label{limsup} With notation as above,
\begin{displaymath}\limsup_{T \rightarrow \infty} \frac{1}{T} \int_0^T X(t) dt \le \lambda
\end{displaymath} with probability 1 for all $\lambda >
E\eta/E\xi$.\end{Cor}

In order to prove this corollary, we need the following technical lemma:

\begin{lemma}\label{tech} Let $0 < \gamma < 1$. Let $U: \R^+ \rightarrow
\R^+$ be such that for
all sequences $\{a_n\}_{n=0}^{\infty}$ with $\sum_{n=0}^{\infty}
\gamma^{a_n}$ convergent, \begin{displaymath} \limsup_{n \rightarrow \infty} U(a_n)
\le c, \end{displaymath} for some $c>0$. Then \begin{displaymath}
\limsup_{T \rightarrow \infty} U(T) \le c. \end{displaymath} \end{lemma}

\noindent
\textbf{Proof:} We proceed by contradiction. Suppose $\limsup_{T
\rightarrow \infty} U(T) > c$. Then, there is a sequence of time $t_n$,
$t_n \rightarrow \infty$, such that $U(t_n) > c$. Take a subsequence
$t_{n_k}$, where $n_k$ is such that $t_n > k$ for all $n \geq n_k$. Such
a subsequence exists since $t_n$ diverges. Now, letting $a_k = t_{n_k}$,
note that $a_k > k$, so $\sum_{k=0}^{\infty} \gamma^{a_k}$ is convergent.
So $\limsup_{k \rightarrow \infty} U(a_k) \le c$. But by definition,
$U(a_k) > c$ for all $k$. This is a contradiction.\qed\medskip

We now proceed with the proof of Corollary~\ref{limsup}. Let $U(T) =
\frac{1}{T}\int_{0}^T X(t)dt$. Let $\gamma$ be as in the conclusion of
Proposition~\ref{gendev}. Then, for any sequence $a_n$ we have
\begin{displaymath} \sum_{n=0}^{\infty} P\left(U(a_n) > \lambda\right)
\le \sum_{n=0}^{\infty} \gamma^{a_n}.
\end{displaymath} Thus, if $\sum_{n=0}^{\infty} \gamma^{a_n}$ converges, by the
Borel-Cantelli lemma, $\limsup_{n \rightarrow \infty} U(a_n) \le \lambda$,
with probability one.
Applying Lemma~\ref{tech}, we have our result. \qed\medskip

\section{Proofs of main results}\label{proofs}

We fix the following notation for the rest of this section: fix a $\delta>0$, and fix $q \in Q$. For $h \in SL(2, \R)$, we define $V_{\delta, q}(h) = V_{\delta}(hq)$. $V_{\delta,q}$ is $K$-invariant, and thus can be viewed as a function on $\h^2$. We define the sets $C_{\delta,l}(q) = \{z \in \h^2: V_q(z) \le l\}$. In the rest of this section we work in this $\h^2$, identifying $g_t r_{\theta} q \in Q$ with $i. g_t r_{\theta} \in \h^2$. For notational convenience, we drop the subscripts $\delta$ and $q$, and write $V$ and $C_l$ for $V_{\delta,q}, C_{\delta,l}(q)$. Furthermore, all distances are measured in $\h^2$.

\vspace{.1in}

\noindent\textbf{Proof of Theorem~\ref{flowreturn}:} Let $q \notin C_l$. We want to estimate the measure of the sets
$$B^{\prime}(T, l,
q) = \{\theta: i.g_t r_{\theta} \notin C_l, 0 \le t \le T\}.$$ For
technical reasons, we will instead study the sets $$B(T):=B(T, l, q) =
\{\theta: \exists \phi \in B^{\prime}(T, l, q) \mbox{ such that }d(i.g_T
r_{\theta}, i.g_T
r_{\phi}) \le 3\tau\},$$ where we will specify $\tau$ shortly. Note that,
by definition, and logsmoothness of $V$, there is some $a_{\tau} \geq 1$ such that $\phi \in B(T, l, q)$ implies that $V(i.g_t
r_{\phi}) > l/a_{\tau} =: l^{\prime\prime}$ for all $0 \le t \le T$.

Let $B_{n\tau}
= B(n\tau, l, q)$, and $p_{n\tau} = \nu(B_{n\tau})$. We have \begin{displaymath}
l^{\prime\prime} p_{n\tau}
\le \int_{B_{n\tau}} V(i.g_{n\tau}r_{\theta}) d\theta  =: D_{n\tau}.
\end{displaymath}

Our main lemma is as follows 

\begin{lemma}\label{average}  For all $\delta >0$, and all $\tau$ sufficiently large there are constants $c = c(\delta), b= b(\tau, \delta)$, so that \begin{displaymath} D_{n\tau} \le ce^{-(1-\delta)\tau} + b\end{displaymath} \end{lemma}

\noindent\textbf{Proof:} Note that by definition $B_{n\tau}$ is a union of arcs of the form $I(n\tau)$, and as such, both Lemma~\ref{polar2} and Corollary~\ref{polar3} apply (we are also using the fact that $V$ is logsmooth). Fix $\sigma>1$, and let $\tau$ be such that we can apply corollary~\ref{polar3}. Setting $c = \sigma c^{\prime}\tilde{c}$ and $b = \sigma c^{\prime} \tilde{b}$, we obtain our result.\qed\medskip

Let $\tau$ be large enough so that we can apply Lemma~\ref{average}. Proceeding as in the proof of Proposition~\ref{walkgen} part (1), we obtain $$p_{n\tau}(q) \le \frac{V_{\delta}(q)}{l^{\prime\prime}} \left(ce^{-(1-\delta)\tau} + \frac{b}{l^{\prime\prime}}\right)^n.$$ Let $\tau_0 >0$ be such that  $ce^{-(1-\delta)\tau_0} <1$. 

Let $t > \tau_0$, and  let $$l_0 \geq a \sup_{\tau_0 \le \tau \le 2\tau_0} \frac{b}{(1- ce^{1-\delta)\tau})}.$$ Let $$a = \sup_{\tau_0 \le \tau \le 2\tau_0} a(\tau),$$ and set $l^{\prime} = l/a$.

Let $$\delta^{\prime} = \delta + \sup_{\tau_0 \le \tau \le 2\tau_0} \frac{1}{\tau}\ln\left(c + \frac{b}{l^{\prime}}e^{(1-\delta)\tau}\right).$$ It is easy to check that $\delta^{\prime} < 1$, and that it is decreasing as a function of $l$. There is some $\tau_0 \le \tau \le 2\tau_0$ and $n \in \N$ such that $t = {n\tau}$. We have $$p_{n\tau} \le \frac{V_{\delta}(q)}{l^{\prime}}\left(ce^{-(1-\delta)\tau} + \frac{b}{l^{\prime}}\right)^n.$$ Rewriting this, we obtain $$p_t \le a\frac{V_{\delta}(q)}{l}e^{-(1-\delta^{\prime})t}.$$\qed\medskip

\vspace{.1in}
\noindent 
\textbf{Proof of Theorem~\ref{sojourn}:} Let $q \in Q$. Consider the circle of radius $S+T$, $\{i.g_{S+T} r_{\theta}: 0 \le \theta <2\pi\}$. Given $T>0$, we want to show that the set $B = B_{S, T} (q,l) =
\{\theta: i.g_t r_{\theta}  \notin C_l, S \le t \le S+T\}$ has
exponentially small measure in $T$ for sufficiently large $l, S,$ and $T$. 

Given $\theta_0 \in B$, let $z_0 = i. g_{S} r_{\theta_0}$, and consider the
circle $\{i.g_T r_{\phi} g_{S}r_{\theta_0} : 0 \le \phi < 2\pi\}$ of radius $T$
around it. By Theorem~\ref{flowreturn}, we know that for
most (the complement is exponentially small in $T$) directions $\phi$ on
this circle, $i.g_t r_{\phi} g_{S} r_{\theta_0} \in C_l$ for some $t < T$. 

Our idea is as follows: there is a small neighborhood $U$ of $\theta_0$
such that each geodesic trajectory $\{i.g_t r_{\theta}\}_{t=0}^{S+T}, \theta
\in U$ is closely shadowed by a piecewise geodesic of the form
$\gamma = \gamma_{\theta_0, \phi}$,  $\theta = \Psi_{S, T}(\phi)$, where
\begin{equation}\gamma(t) = \left\{\begin{array}{ll} i.g_t r_{\theta_0}
 & 0 \le t \le S \\ i.g_{t-S}r_{\phi}g_{S}r_{\theta_0} & S < t <
S+T\end{array}\right.\nonumber\end{equation}

By closely shadowed, we mean that $d(\gamma(t), i.g_t r_{\theta})$ is small for all $t$. Thus, if $V(\gamma(t)) \le l$ for some $S \le S+T$, we have that $V(i.g_t r_{\theta}) \le \tilde{l}$, for some $\tilde{l} > l$.

Now, for all but a small set of $\phi$, $\gamma(t) \in C_l$ for some $S
\le t \le S+T$. Thus, in a small neighborhood of $\theta_0$, we have a
(large-proportioned) collection of angles which are not in $\tilde{B} = B_{S,T}(q, \tilde{l})$. 

To make this rigorous, fix $\epsilon, \delta >0$, and let $A: = A(\theta_0) = B_{0,T}(z_0, l) = \{\phi: V(i.g_t r_{\phi}g_{S}r_{\theta_0})>l, 0 \le t \le T\}$. Let $d$ be the maximum thickness of a hyperbolic triangle. By logsmoothness of $V$, there is a $b \geq 1$ such that $d(z_1, z_2) \le d$ implies $\frac{V(z_1)}{V(z_2)} \le b$. Let $l_1 = bl_0$, where $l_0$ is as in the conclusion of Theorem~\ref{flowreturn}. For any $l$, let $\hat{l} = l/b$.

Setting $\delta^{\prime\prime} =\delta^{\prime}(\hat{l}, \delta)$ and $l^{\prime} = \hat{l}/a = l/ba$, where $a = a(\delta)$ is as in Theorem~\ref{flowreturn} yields $$\nu(A) \le \frac{V(z_0)}{l^{\prime}}e^{-(1-\delta^{\prime\prime})t}.$$ Let $T, S$ be large enough so that we can apply Lemma~\ref{shadow} with $\epsilon$, and $t_2= T$, $t_1 =
S$. This gives a neighborhood $U$ of $\theta_0$ in which the proportion
of angles $$\frac{\nu(\Psi(A) \cap U)}{\nu(U)} \le 4(1+\epsilon)\sup_{\theta \in [0, 2\pi)}\frac{V(i.g_S r_{\theta})}{l^{\prime}}e^{-(1-\delta^{\prime\prime})t}.$$ Now, by the thinness of triangles in hyperbolic geometry,
given $\theta \in U$, there is a $\phi$ with $\Psi(\phi) = \theta$ such that
 \begin{displaymath} d(i.g_t r_{\theta},
\gamma_{\theta_0, \phi}(t)) \le d \end{displaymath} for $t \in [0,
S+T]$.

For all $\phi \notin A$, we have \begin{displaymath}\gamma_{\theta_0,
\phi} (t_0) \in
C_{\hat{l}},\end{displaymath} for some $S+T > t_0 > S$. Thus,
\begin{displaymath} g_{t_0} r_{\theta} q \in C_l,\end{displaymath} so
$\theta \notin B$. 

Thus, given any $\theta_0 \in B$, we have produced a neighborhood
$U$ s.t. \begin{displaymath} \frac{\nu(B^C \cap U)}{\nu(U)} > 1- 4(1+\epsilon)\sup_{\theta \in [0, 2\pi)}\frac{V(i.g_S r_{\theta})}{l^{\prime}}e^{-(1-\delta^{\prime\prime})t}. \end{displaymath}

To complete the proof, we need the following standard lemma (see, for example,~\cite{Evans}):

\begin{lemma}\label{besic} Let $B \subset [0,2\pi)$ be a measurable set such
that for all $b \in B$, there is a $\delta_b >0$, so that $U_b =
[b-\delta_b, b+\delta_b] \subset [0,2\pi)$ satisfies \begin{displaymath} \frac{\nu(U_b
\cap B)}{\nu(U_b)} < \epsilon.\end{displaymath} Then
\begin{displaymath} \nu(B) \le 2 \epsilon.\end{displaymath}\end{lemma}
\noindent

\noindent Applying the lemma to our set $B$, we obtain $$\nu(B) \le
8(1+\epsilon)\sup_{\theta \in [0, 2\pi)}\frac{V(i.g_S r_{\theta})}{l^{\prime}}e^{-(1-\delta^{\prime\prime})t}.$$\qed\medskip

\vspace{.1in}

\noindent\textbf{Proof of Theorem~\ref{proportion}:} Our
strategy is as follows: Given a direction
$\theta$, consider the succesive departures and returns of the geodesic
trajectory $\{i.g_t r_{\theta}\}_{t \geq 0}$ to the compact set $C_l$.
Theorem~\ref{flowreturn} implies that the probability any departure is
long is small, and thus, we can try and apply Proposition~\ref{gendev}, to
the `random variables' given by the length of sojourns inside and outside
the compact set. 

We proceed as follows: Let $d$ be as in the proof of
Theorem~\ref{sojourn}. Fix $\delta >0$, and let $l_0$ be as in Theorem~\ref{flowreturn}. Let $l > l_0$ be such that  $d(C_l^{c},C_{l_0}) > 2d$, and define $C = C(l) = d(C_l^{c},C_{l_0}) - 2d$.  Define
$t_0(\theta) =0$, and set \begin{displaymath} t_{2n}(\theta) = \inf \{t>
t_{2n-1}: \exists \phi \mbox{ such that }d(i.g_t r_{\phi}, i.g_t
r_{\theta}) < d, i.g_t r_{\phi} \in
C_{l_0}\}\end{displaymath} and
\begin{displaymath} t_{2n+1}(\theta) = \inf \{t> t_{2n}: \exists \phi
\mbox{ such that } d (i.g_t r_{\phi}, i.g_t r_{\theta}) < d,
i.g_t
r_{\phi} \notin C_l\},\end{displaymath} for $n \geq 0$.
Define
$\tau_i(\theta) = t_i - t_{i-1}$. Now fix $C^{\prime}>0$, and define auxiliary functions $\tau^{\prime}_i$ by 
$$\tau^{\prime}_{2i} := \left\{\begin{array}{ll} 0
 & \tau_{2i} \le C^{\prime} \\ \tau_{2i} & \tau_{2i} >C^{\prime}\end{array}\right.$$
and
$$\tau^{\prime}_{2i-1} := \left\{\begin{array}{ll} \tau_{2i-1}+\tau_{2i}
 & \tau_{2i} \le C^{\prime} \\ \tau_{2i-1} & \tau_{2i} >C^{\prime}\end{array}\right.$$
 
Define $t^{\prime}_i := \sum_{j=1}^i \tau^{\prime}_j$. Note that $\tau_{2i} + \tau_{2i-1} = \tau^{\prime}_{2i} + \tau^{\prime}_{2i-1}$, so $t^{\prime}_{2i} = t_{2i}$.

\begin{equation}X(t) = \left\{
\begin{array}{ll}  1 & t^{\prime}_{2i-1} \le t
< t^{\prime}_{2i}\\ 0 & \mbox{otherwise}\end{array}\right.\nonumber\end{equation}

The $t_i$'s should be thought of as the entry and departure times of the trajectory $\{i.g_t r_{\theta}\}_{t \geq 0}$ from our compact sets. For technical reasons, they are defined as the first time when a \emph{nearby} trajectory leaves a larger compact set  ($C_l$), or re-enters a smaller one ($C_{l_0}$). The $\tau_{2i}$'s measure time spent after departing the larger set before returning to the smaller, and the $\tau_{2i-1}$'s measure the time spent after returning to the smaller before departing the larger. The auxiliary $\tau^{\prime}_i$'s are defined to exclude short ($\le C^{\prime}$) sojourns. Thus, if $t^{\prime}_{2i} < t < t^{\prime}_{2i+1}$ we are within distance $C^{\prime}$ of the larger compact set, and thus still within a compact set. Thus, for all such $t$, $V(i.g_t r_{\theta}) \le \tilde{l} = \tilde{l}(C^{\prime})$, i.e., \begin{displaymath}
X(t) \geq \chi_{C_{\tilde{l}}^{c}}(i.g_t r_{\theta} ).\end{displaymath} Thus, it suffices to show that for all $\lambda >0$ there is a $C^{\prime} >0$ so that
$\nu\left(\frac{1}{T} \int_{0}^{T} X(t)dt > \lambda\right)$ decays exponentially in
$T$. 

We will apply Proposition~\ref{gendev} to $\{\tau^{\prime}_i\}_{i \geq 0}$, with $\Omega = S^1$, $\F$ the
standard
$\sigma$-algebra, and $P=\nu$ the Haar measure. Let $\F_n = \sigma(t^{\prime}_1,
\ldots, t^{\prime}_n)$ be the $\sigma$-algebra generated by $t^{\prime}_1, \ldots, t^{\prime}_n$.
Clearly $\tau^{\prime}_n$ is $\F_n$-measurable. Note that
since $d(C_l^{c}, C_{l_0}) = 2d + C$, \begin{displaymath}
\tau^{\prime}_{2i+1} \geq \tau_{2i+1}
= t_{2i+1} - t_{2i} > C \end{displaymath} (since $i.g_{t_{2i}}r_{\theta}$
is within distance $d$ of $C_{l_0}$ and
$i.g_{t_{2i+1}}r_{\theta}$
is within distance $d$ of $C_{l}^c$). Thus, condition 2 of the
proposition
is satisfied, with $\xi = C$

It remains to check condition 1. We will show
\begin{equation}\label{decay} \nu(\tau^{\prime}_{2i} > t | \F_{2i-1}) \le a_1
e^{-a_2t}\end{equation} for some $a_1, a_2 >0$ and all $t>C^{\prime}$.

By definition, $t_i^{\prime-1} (x)$ is an interval, and thus $\F_i$ is
generated by these intervals. Fix $\theta$. Let \begin{displaymath}
I_n(\theta) = \{ \phi: t_i^{\prime}(\phi) = t_i^{\prime}(\theta) \mbox{ for all } 0 \le i \le
n\}.\end{displaymath}

It suffices to show that there are $a_1, a_2 >0$ such that, for all $t>C^{\prime}$,
\begin{equation}\label{decay2} \nu\left(\tau^{\prime}_{2n} > t|
I_{2n-1}(\theta)\right) \le a_1
e^{-a_2t}.\end{equation} Fix $\phi \in I_{2n-1}(\theta)$ with
$\tau_{2i}(\phi) > t$. Let $z_{\phi} = i.g_{t_{2i-1}} r_{\phi}$. Note that
this point is within distance $d$ of the boundary of $C_l^c$, thus it
is both outside $C_{l_0}$ and still contained within a compact set.
Consider the circle of radius $t$ around $z_{\phi}$ and the associated map
$\Psi = \Psi_{t, t_{2i-1}}$ back to the circle at $i$.

Let $A = A(\phi) = \{\theta: i.g_s r_{\theta}g_{t_{2i-1}}r_{\phi} \notin
C_{l_0},  0 \le s \le t \}$. Since $z_{\phi}$ is contained in a
compact set, we can pick $c_1 =c_1(l), c_2 =c_2(l_0)$ independent of $\phi$ such that
\begin{displaymath} \nu(A) \le c_1e^{-c_2t},\end{displaymath} for all $t$ sufficiently large. By picking $C^{\prime}$ large enough, we can get this to hold for all $t > C^{\prime}$. For the rest of this section, let $t>C^{\prime}$.

Applying lemma~\ref{shadow} with $\epsilon =1$, we obtain a neighborhood $U$ of $\phi$, with
\begin{displaymath}\frac{\nu(\Psi A \cap U)}{\nu(U)} \le 16c_1
e^{-c_2t}.\end{displaymath}

Note that $U \subset \chi_{2n-1}(\theta)$, since $\forall 0 \le t \le
t_{2i-1}, \theta^{\prime} \in U$, \begin{displaymath} d(i.g_t
r_{\theta^{\prime}}, i.g_t r_{\theta}) < d.\end{displaymath}

Finally, observe that $\tau_{2i}(\theta^{\prime}) < t$ for all $\theta
^{\prime} \in U- \Psi(A(\phi))$, since $\theta^{\prime} \notin \Psi(A)$ implies that there is some $0 \le s \le t$ such that $V(i.g_s r_{\psi}g_{t_{2i-1}}r_{\phi}) \le l_0$, with $\Psi (\psi) = \theta^{\prime}$. Thus, since $d(i.g_s r_{\psi}g_{t_{2i-1}}r_{\phi}, i.g_{s+t_{2i-1}}r_{\theta^{\prime}}) \le d$, we have $V(i.g_{s+t_{2i-1}}r_{\theta^{\prime}}) \le l$, and thus $\tau_{2i}(\theta^{\prime}) < t$. Once again
applying Lemma~\ref{besic}, we obtain equation~\eqref{decay2}, with $a_1 = 32c_1$ and $a_2 = c_2$. 

Let $C^{\prime}$ be large enough so that $a_1e^{-a_2C^{\prime}} <1$. Let $\eta$ be a non-negative function on $S^1$ such that $$\nu\{\theta: \eta(\theta) = 0\} = a_1e^{-a_2C^{\prime}}$$ and $$\nu\{\theta: \eta(\theta) > t\} = a_1e^{-a_2t}$$ for all $t > C^{\prime}$. Condition 1 is then clearly satisfied, since $\nu(\tau_{2i}^{\prime} =0) \geq a_1e^{-a_2C^{\prime}}$, and $\nu(\tau^{\prime}_{2i} > t) \le a_1e^{-a_2 t}$ for $t> C^{\prime}$. Now, note that $E(\eta) = \frac{a_1}{a_2}e^{-a_2 C^{\prime}}$, and $E(\xi) = C$. Thus, by enlarging $C^{\prime}$, we can make $E(\eta)/E(\xi)$ arbitrarily small. Precisely for any $\lambda >0$, take $C^{\prime}$ so that $$E(\eta)/E(\xi) \le \lambda.$$ Then, setting $\tilde{l} = \tilde{l}(C^{\prime})$ and applying Proposition~\ref{gendev}, we obtain that there is a $\gamma <1$ so that $$\nu\{\theta: \frac{1}{T}|\{0 \le t \le T:g_t r_{\theta} q \notin C_{\tilde{l}}\}|
> \lambda\} \le \nu\left(\frac{1}{T} \int_{0}^{T} X(t)dt > \lambda\right) \le \gamma^T,$$ for all $T$ sufficiently large.  \qed\medskip

\noindent
\textbf{Proof of Corollary~\ref{for}:} Applying Corollary~\ref{limsup} to
$X(t)$, we obtain our result.\qed \medskip

\noindent
\textbf{Acknowledgements:} I would like to thank my advisor, Professor Alex Eskin, for his guidance
throughout this project. I would also like to thank Professors Howard
Masur, Giovanni Forni, Steven Lalley and Krishna Athreya for valuable discussions, and my colleague Matthew Day for help with hyperbolic trigonometry. Thanks are also due to the anonymous referee, whose detailed comments and suggestions greatly improved this paper.

\end{document}